\documentclass[reqno]{amsart}
\usepackage{setspace}
\usepackage{graphicx}
\usepackage{amsmath}
\usepackage{amsfonts}
\usepackage{amssymb}
\begin{document}

\centerline{\Huge \bf On Tomaszewski's Cube Vertices Problem}

\bigskip
\centerline{\Large Yiming Li, Yuqin Zhang and Miao Fu\footnote{Corresponding author. }}

{\large
\vspace{0.8cm}
\centerline{\begin{minipage}{12.8cm}
{\bf Abstract.} The following assertion was equivalent to a conjecture proposed by B. Tomaszewski :  {\it Let $C$ be an $n$-dimensional unit cube and let $H$ be a plank of thickness $1$, both are centered at the origin,  then no matter how to turn the cube around, $C\cap H$ contains at least half of the cube's $2^n$ vertices.} A lower bound for the number of the vertices of $C$ in $C\cap H$ was obtained.
\end{minipage}}

\vspace{0.6cm}

\vspace{1cm}
\noindent
{\LARGE\bf 1. Introduction}

\bigskip
\noindent

As a fundamental object in geometry, in combinatorics, in binary codes, and in optimization, the $n$-dimensional unit cube itself is still far from being well understood. Many basic problems about it are still open (see \cite{Zong-05, Zong-06}). For example, what is the maximum area of its cross sections? What is the maximum area of its projections? What is the maximum volume of a simplex inscribed in it? What is the smallest number of simplices to triangulate it?

For convenience, we write
$$B_n=\left\{(x_1, x_2, \ldots , x_n):\ x_1^2+x_2^2+\ldots +x_n^2\le 1 \right\},$$
$$C_n=\left\{(x_1, x_2, \ldots , x_n):\ \left|x_i\right|\le 1, \ i=1, 2, \ldots , n \right\},$$
and
$$M_n=\left\{(x_1, x_2, \ldots , x_n):\ x_i=\pm 1, \ i=1, 2, \ldots , n \right\}.$$
In other words, $B_n$ is an $n$-dimensional unit ball centered at the origin, $C_n$ is an $n$-dimensional cube of edge length $2$ centered at the origin, and $M_n$ is the set of the vertices of $C_n$. Furthermore, as usual, let $int(K)$, $\partial (K)$ and $conv(K)$ denote the interior, the boundary and the convex hull of a set $K$, respectively, and let $\sharp\{X\}$ denote the number of the points in a discrete set $X$, let $\langle {\bf x}, {\bf y}\rangle$ denote the inner product of ${\bf x}$ and ${\bf y}$.

In 1986, B. Tomaszewski proposed the following problem, which presented in \cite{Guy-86}: {\it Considers $n$ real numbers $a_1, \ldots, a_n$, such that $\sum_{i=1}^{n}{a_{i}}^{2}=1$. Of the $2^n$ expressions $|{\epsilon_1}{a_1}+ \ldots +{\epsilon_n}{a_n}|$ with $\epsilon_i=1$ $or -1$, $1\leq i \leq n$, can there be more with value $>1$ than with value $\leq 1$ }?

\medskip
By using the language of probability theory, R. Holzman and D. J. Kleitman \cite{Holzman-92} reformulated this problem as follows: {\it Let $X=\sum_{i=1}^{n}a_{i}x_{i}$, where $\{x_{i}\}$ are uniformly distributed in $\{-1, 1\}$ and independent. Is it true that
$$Pr\{|X|\leq 1\} \geq {{1}\over{2}}\eqno(1)$$ holds for arbitrary $n$ real numbers $a_1, a_2, \ldots , a_n$ which satisfied $\sum_{i=1}^{n}{a_{i}}^{2}=1$ }? And then he proved that $$Pr\{|X|\leq 1\} \geq {{3}\over{8}}$$  instead of (1). Recently, N. Keller and O. Klein \cite{Keller-21} proved proposition (1) by a much more complicated analysis, thus Tomaszewski's problem is solved. For other details in probability theory perspective we refer to \cite{Boppana-17}, \cite{Boppana-20}.

\medskip
 However, Tomaszewski's problem also can be reformulated in combinatorial and convex geometry means: {\it For every ${\bf u}\in \partial (B_n)$, the inequality
$$|\langle {\bf x}, {\bf u}\rangle |\le 1$$
has at least $2^{n-1}$ solutions in $M_n$.} As a equivalent form of (1), there are no research try to attack problem in this perspective. We present a combinatorial and convex geometry method to prove the following weaker result:

\medskip\noindent
{\bf Theorem 1.} {\it For every $n$-dimensional plank $P_n$ bounded by two parallel hyperplanes, both are tangent to $B_n$, we have}
$$\sharp \{P_n\cap M_n \}\ge {{2^{n-1}}\over {\sqrt{n}}}.$$

\vspace{0.6cm}
\noindent
{\LARGE\bf 2. Proof of Theorem 1}

\bigskip
\noindent
First, we prove the following lemma.

\medskip\noindent
{\bf Lemma 1.} {\it Let $L_n$ be an half space bounded by an hyperplane $T_n$. If $int (B_n)\cap L_n=\emptyset$, then}
$$\sharp \{int (L_n)\cap M_n \}< \left(2-1/\sqrt{n}\right)2^{n-2}.$$

\medskip\noindent
{\bf Proof.} Without loss of generality, we assume that $T_n$ is tangent to $B_n$ at ${\bf m}_0$. Then, we have
$$L_n\cap B_n=T_n\cap B_n=\{{\bf m}_0\}.$$
Suppose that
$$L_n\cap M_n=\{ {\bf m}_1, {\bf m}_2, \ldots, {\bf m}_k\},$$
clearly we have
$$conv(L_n\cap M_n)\subset L_n.$$
Consequently, we have
$$conv \{ {\bf m}_1, {\bf m}_2, \ldots, {\bf m}_k\}\cap (B_n\setminus \{{\bf m}_0\})=\emptyset .\eqno(2)$$

Then we make two observations:

\smallskip
\noindent
{\bf Observation 1.} If both ${\bf m}_i$ and $-{\bf m}_i$ belong to $L_n\cap M_n$, then by convexity we have
$${\bf o}={1\over 2}({\bf m}_i+(-{\bf m}_i))\in conv(L_n\cap M_n)\subset L_n,$$
which contradicts to (2).

\smallskip
\noindent
{\bf Observation 2.} If both ${\bf m}_i=(a_1,\ldots, a_{j-1}, a_j, a_{j+1}, \ldots, a_n)$ and ${\bf m}'_i=(-a_1,\ldots, -a_{j-1},$ $a_j, -a_{j+1}, \ldots, -a_n)$ belong to $L_n\cap M_n,$ then by convexity we have
$$(0, \ldots, 0, a_j, 0, \ldots , 0)={1\over 2}({\bf m}_i+{\bf m}'_i)\in conv(L_n\cap M_n)\subset L_n.$$
Since $a_j=\pm 1$, $T_n$ must be tangent to $B_n$ at ${\bf m}_{0}=\left(0,\ldots,0,a_{j},0,\ldots,0 \right)$. Then, we get
$$\sharp \{int (L_n)\cap M_n \}=0.$$

Let us consider the first coordinate of the vertices of $C_n$ in $L_n\cap M_n=\{ {\bf m}_1, {\bf m}_2, \ldots, {\bf m}_k\}$. Suppose $m$ of them are $1$ and therefore $k-m$ of them are $-1$. For convenience, we write
$$M'_n=\{ (x_1, x_2, \ldots , x_n)\in M_n: \ x_1=1\}$$
and
$$M^*_n=\{ (x_1, x_2, \ldots , x_n)\in L_n\cap M_n:\ x_1=1\}.$$
Then we define a map $\pi $ from $M^*_n$ to $M'_n$ by
$$(1, x_2, \ldots , x_n) \longmapsto (1, -x_2, \ldots , -x_n).$$
It is easy to see that the map is one-to-one and
$$\sharp \{M^*_n\cup \pi (M^*_n)\}\le \sharp \{M'_n\}=2^{n-1}.$$

Now, we prove the lemma by considering two cases:

\medskip\noindent
{\bf Case 1.} $\sharp \{M^*_n\}=m>2^{n-2}.$

Then, we have
$$\left\{\begin{array}{ll}
2^{n-1}\hspace{-0.2cm} &< \sharp \{\pi (M^*_n)\} +\sharp \{M^*_n\},\\
2^{n-1}\hspace{-0.2cm} &\ge \sharp \{\pi (M^*_n) \cup M^*_n\}\\
\end{array}\right.$$
and therefore, by De Morgan's laws,
$$M^*_n\cap \pi (M^*_n)\not=\emptyset.$$
In other words, $M^*_n$ has two points satisfying the condition of Observation 2. Then, we get
$$\sharp \{int (L_n)\cap M_n \}=0.$$

\noindent
{\bf Case 2.} $\sharp \{M^*_n\}=m\le 2^{n-2}.$ {\it Then, without loss of generality, we may assume that
$$k-m\le 2^{n-2}\eqno(3)$$
holds as well.}

If, on the contrary to the lemma,
$$k\ge \left(2-1/\sqrt{n}\right)2^{n-2},\eqno(4)$$
then one can deduce from (3) and (4) that
$$2^{n-2}\ge m\ge k-2^{n-2}\ge \left(1-1/\sqrt{n}\right)2^{n-2}.\eqno(5)$$
Similarly, we also have
$$2^{n-2}\ge k-m\ge \left(1-1/\sqrt{n}\right)2^{n-2}.\eqno(6)$$

We write that
$${\bf w}=(w_1, w_2, \ldots, w_n)={\bf m}_1+{\bf m}_2+\ldots +{\bf m}_k.$$
Then, by (5) and (6), we get
$$-{{2^{n-2}}\over {\sqrt{n}}}\le w_1=2m-k\le {{2^{n-2}}\over {\sqrt{n}}}.$$
Similarly, one can deduce that
$$-{{2^{n-2}}\over {\sqrt{n}}}\le w_i\le {{2^{n-2}}\over {\sqrt{n}}}\eqno(7)$$
holds for all coordinates $w_i$ of ${\bf w}$. Consequently, by (7) and (4) we get
$$\|{\bf w}\|^2=\sum_{i=1}^nw_i^2\le n\left({{2^{n-2}}\over {\sqrt{n}}}\right)^2 =2^{2(n-2)},$$
$$\left\|\mbox{${1\over k}$}{\bf w}\right\|^2\le {{2^{2(n-2)}}\over {k^2}}\le \left(2-1/\sqrt{n}\right)^{-2}<1$$
and therefore
$$\mbox{${1\over k}$}{\bf w}\in int (B_n),$$
which contradicts to (2). Thus, in this case, we must have
$$\sharp \{int (L_n)\cap M_n \} \le k< \left(2-1/\sqrt{n}\right)2^{n-2}.$$

As a conclusion of the two cases, the lemma is proved. \hfill{$\Box$}

\medskip
\noindent
{\bf Proof of Theorem 1.} Clearly, by Lemma 1 we have
$$\sharp \{P_n\cap M_n \}= 2^n-2\cdot \sharp \{int (L_n)\cap M_n \} \ge {{2^{n-1}}\over {\sqrt{n}}}.$$
Theorem 1 is proved. \hfill{$\Box$}

\medskip
\noindent
{\bf Remark.} Although Theorem $1$ is weaker than the result of N. Keller and O. Klein, we believed that by more deeper observation and analysis, Tomaszewski's problem can also solved by combinatorial and convex geometry method independently as probability theory.

\vspace{0.6cm}\noindent
{\bf Acknowledgements.} This work is supported by the National Natural Science Foundation of China (NSFC11921001) and the National Key Research and Development Program of China (2018YFA0704701) and the Scholarship Council of China. The author is grateful to professor C. Zong for his valuable supervision and discussion.

\bibliographystyle{amsplain}

\bigskip\medskip\noindent
Yiming Li

\noindent
Center for Applied Mathematics\\
Tianjin University\\
Tianjin 300072\\
P. R. China

\noindent
Email: xiaozhuang@tju.edu.cn

\bigskip\medskip\noindent
Yuqin Zhang

\noindent
School of Mathematics\\
Tianjin University\\
Tianjin 300072\\
P. R. China

\noindent
Email: yuqinzhang@tju.edu.cn

\bigskip\medskip\noindent
Miao Fu

\noindent
Center for Applied Mathematics \\
Tianjin University\\
Tianjin, 300072\\
P. R. China

\noindent
Email: miaofu@tju.edu.cn
\end{document}